\documentclass[11pt]{article}
\usepackage{amsmath,amsfonts,amssymb,amsthm,amscd}
\usepackage{color}
\numberwithin{equation}{section} \allowdisplaybreaks
 \usepackage[all]{xy}

\newtheorem{theorem}{\sc Theorem}[section]
\newtheorem{lemma}[theorem]{\sc Lemma}
\newtheorem{proposition}[theorem]{\sc Proposition}
\newtheorem{corollary}[theorem]{\sc Corollary}
\newtheorem{definition}[theorem]{\sc Definition}
\newtheorem{example}[theorem]{\sc Example}
\newtheorem{remark}[theorem]{\sc Remark}

\newcommand{\bet}{\begin{theorem}}
\newcommand{\eet}{\end{theorem}}
\newcommand{\blm}{\begin{lemma}}
\newcommand{\elm}{\end{lemma}}
\newcommand{\bprop}{\begin{proposition}}
\newcommand{\eprop}{\end{proposition}}
\newcommand{\bcor}{\begin{corollary}}
\newcommand{\ecor}{\end{corollary}}
\newcommand{\bdf}{\begin{definition}\rm}
\newcommand{\edf}{\end{definition}}
\newcommand{\bp}{\begin{proof}}
\newcommand{\ep}{\end{proof}}
\newcommand{\bex}{\begin{example}\rm}
\newcommand{\eex}{\end{example}}
\newcommand{\bremark}{\begin{remark}\rm}
\newcommand{\eremark}{\end{remark}}

\newcommand{\norm}[1]{ \| #1 \| }

\newtheorem{df}{Definition}[section]

\newtheorem{lem}[df] {Lemma}

\begin{document}

\title { Index of continuous semi-Fredholm families, regularities  and  semiregularities }

\author{ M. Berkani}

\date{}

\maketitle

\begin{abstract}

In this paper, we define and index for   continuous families of
semi-Fredholm bounded liner operators. Moreover, we study various
regularities and semiregularities of continuous families of
bounded linear operators.

\end{abstract}
\footnotetext{\hspace{-7pt}2010 {\em Mathematics Subject
Classification\/}:  47A53, 47A60 \baselineskip=18pt\newline\indent
 }

\section{Introduction}
In this part,  we will introduce some definitions  and notations
that will be needed in the next sections.

In \cite{MM}, the concept  of regularity  was defined by the following :\\
\vspace{-5 mm }

\bdf \label{REG}
 Let $A$ be a unital Banach algebra, with unit $\bf{e}$.
 A non-empty subset $\mathbf{R}$ of $A$ is called a
regularity if it satisfies the following conditions:\\
 \noindent i) If $ a \in A$ and  $n\geq 1,$ then  $
a \in {\bf R} $ if and only if $a^n \in {\bf R}. $\\
 \noindent ii) If $ a,b,c,d \in L(X) $  are mutually commuting
 elements of $A$
satisfying  $ ac+bd= \bf{e}, $ then $ab \in {\bf R} $  if and only
if
 $a,b \in {\bf R} $ .
\edf

\noindent  A regularity ${\bf R} $ defines in a natural way a
spectrum by $\sigma_{\bf R} (a)  = \{ \lambda \in \mathbb{C}:
a-\lambda \bf{e}  \notin{\bf R} \},$  for  $a \in A.$
 Moreover the spectrum
$\sigma_{\bf R}$ satisfies the spectral mapping
theorem.\\
\noindent Let $X$ be an infinite dimensional Banach space and
$L(X)$ be the Banach algebra of bounded linear operators acting
 $X.$ Various regularities
of the Banach algebra $L(X)$ were studied in \cite{MM} and
\cite{P8}. Here, we will consider similar regularities but for
continuous families  of bounded linear operators.

 \noindent In the sequel if $E$  and $F$ are two
vector spaces, the notation   $E  \simeq   F $  will  mean that
$E$ and  $F$ are isomorphic. If $E$ and $F$ are vector subspaces
of the same vector space $X,$   we will say that $E$ is
essentially included in $F$ and write
 $E  \subset_{e}   F ,$ if there exist a finite
dimensional vector
  subspace  $G$ of $X$ such that
$ E\subset F+G.$ We will say that $E$ is essentially equal to $F$
and write $E =_{e} F $ if
 $E \subset_{e} F $ and  $F  \subset_{e} E $.
 We define also the infimum of the empty set  to be  $ \infty  $.\\
For $ T \in L(X),$ we will denote by $N(T)$ the null space of $T$,
by $ \alpha(T)$ the nullity of $T$,  by $ R(T)$ the range of $T$
and by $\beta(T) $ its defect. If the range $ R(T)$ of $T$ is
closed and $ \alpha(T) < \infty $ (resp. $\beta(T)< \infty $ ),
then $T$ is called an upper semi-Fredholm (resp. a lower
semi-Fredholm) operator. A semi-Fredholm operator is an upper or a
lower semi-Fredholm operator. In the sequel $\Phi_{+}(X)$ (resp.
$\Phi_{-}(X)$) will denote the set of upper (resp. lower)
semi-Fredholm operators and {\bf \scriptsize
S}$\Phi(X)=\Phi_{+}(X) \cup \Phi_{-}(X)$ will denote the set of
semi-Fredholm operators. If both of $ \alpha(T) $ and $\beta(T) $
are finite, then $T $ is called a Fredholm operator and the index
of $T$  is defined by $ ind(T) = \alpha(T) - \beta(T); $ $\Phi(X)$
will denote the set of Fredholm operators. It is well known that
$\Phi_{+}(X), \Phi_{-}(X),$ {\bf \scriptsize S}$\Phi(X)$ and
$\Phi(X)$ are all open subset of $ L(X),$ see for example
\cite[Chapter 4] {CPY}.

\bdf   Let    $ T \in
  L(X), \,\, n\in \mathbb{N} $ and let $c_{n}(T)=dim \, R(T^n)/R(T^{n+1})$.
 The  descent of $T$ is defined by
 $ \delta(T)=inf \{ n:c_{n}(T)=0 \}$=inf $ \{n: R(T^n) = R(T^{n+1})\}$
and the essential descent of $T$ is defined by
 $ \delta_{e}(T)=inf \{ n:c_{n}(T)<\infty \}$= \\ \noindent inf $\{n: R(T^n) =_{e}
 R(T^{n+1})\}.$
\edf

\bdf Let  $ T \in
  L(X),  n\in \mathbb{N}$ and let $c'_{n}(T)=dim \, N(T^{n+1})/N(T^n)$.
Then the  ascent $ a(T) $ of $ T $ is defined by
 $ a(T)=inf \{ n:c'_{n}(T)=0 \}$= inf $\{n : N(T^{n})=N(T^{n+1})\}$
 and the essential ascent of $T$ is defined by
$ a_{e}(T)$= inf$ \{ n:c'_{n}(T)<\infty \}$=\\ \noindent
 inf $\{n: N(T^{n+1})=_{e} N(T^{n})\}.$\\
\edf
 \bdf \cite{GR} {\bf :}
 Let $ T \in L(X)$. The sequence
 $(k_{n}(T))_n$ is defined by
 $ k_{n}(T)=dim[(R(T^n)\cap N(T))/(R(T^{n+1})\cap N(T)),$  for all
 $ n \in \mathbb{N}.$
\edf

\bdf \cite {GR} Let  $ T \in L(X) $ and let $  d \in  \mathbb{N}
$. Then $T$ has a uniform descent for $ n \geq d$ if $R(T) +
N(T^{n}) = R(T) + N(T^{d}) $ for all $ n \geq d.$ If in addition
$R(T) + N(T^{d}) $is closed, then  $T$ is said to have a
topological uniform descent for $ n \geq d$. \edf


\noindent Let $\mathbb{X}$ be a compact topological space and let
$\mathcal{C}( \mathbb{X}, L(X))$ be the algebra of continuous maps
from the topological space $\mathbb{X}$ into the Banach algebra $
L(X).$ If $\mathcal{T}: \mathbb{X} \rightarrow L(X),$ we denote by
$ \mathcal{T}_x$ the image $\mathcal{T}(x)$ of an element  $ x \in
\mathbb{X}.$

\noindent Recall that $\mathcal{C}( \mathbb{X}, L(X))$ is a unital
algebra with the usual properties of addition, scalar
multiplication and multiplication defined by:

$ (\lambda \mathcal{S+ T})_x= \lambda \mathcal{S}_x +
\mathcal{T}_x,   (\mathcal{ST})_x=
\mathcal{S}_x\mathcal{\mathcal{T}}_x,
 \forall \mathcal{(S, T)} \in \mathcal{C}( \mathbb{X}, L(X))^2,  \forall
x \in \mathbb{X},  \forall  \lambda \in \mathbb{C}.$

\noindent The unit element  $ \mathcal{I}$  of the algebra
$\mathcal{C}( \mathbb{X}, L(X))$ is the map defined by:  $
\mathcal{I}_x= I,$  for all $  x \in \mathbb{X},$ where $I$ is the
identity operator. Moreover as $\mathbb{X}$ is compact, if we set
$ \norm{\mathcal{T}}= sup_{x\in \mathbb{X}} \norm {\mathcal{T}_x},
\forall   \, \mathcal{T} \in \mathcal{C}( \mathbb{X}, L(X)),$ then
$\mathcal{C}( \mathbb{X}, L(X))$ equipped with this norm is a
unital Banach algebra.\\

\vspace{2mm}

\section { Index  of continuous semi-Fredholm families }
In \cite{P48}  and \cite{P49} respectively, we defined an
analytical index for continuous Fredholm families of bounded
linear operators
 into  a Hilbert
space (respectively into a Banach space). Here, in the same way as
in \cite{P48}  and \cite{P49},  we define an analytical index for
continuous semi-Fredholm families into a Banach space.

\noindent Let $ \mathbb{X}$ be a topological space. Define an
equivalence relation $ \textbf{\textit{c}} $ on the space
$\mathbb{X}$ by setting that $ x \,\textbf{\textit{c}}\, y,$ if
and only if $ x$ and $y$ belongs to the same connected component
of $\mathbb{X}.$ Let $\mathbb{X}^{\textbf{\textit{c}}} $ be the
quotient space associated to this equivalence relation. We extend
here  the definition \cite [ Definition 1.1]{P49} of the index
of a continuous family of Fredholm operators to the case   of a
continuous family of semi-Fredholm operators.

\bdf \label{index}  The index   of a continuous  family of
semi-Fredholm operators $ \mathcal{T} : \mathbb{X} \rightarrow
L(X),$
   parameterized by a topological space $\mathbb{X}$    is defined by $ ind(\mathcal{T})=
( ind(\mathcal{T}_x))_{\overline{x} \in \mathbb{X}^{c}}.$

 \edf

\noindent Thus the index of  a continuous family of semi-Fredholm
operators $ \mathcal{T}$ is a sequence of elements of  $
\overline{\mathbb{Z}}=\mathbb{Z}\cup \{-\infty, +\infty\},$ which
may be a finite sequence or infinite sequence, depending on the
cardinal $n_c$ of the connected components of $\mathbb{X},$ which
is assumed here to have  at most a countable set of connected
components.

 \bet The index of a continuous family of  semi-Fredholm  operators $\mathcal{T}$  parameterized
by a topological space $\mathbb{X} $ is well defined as an element
of $ \overline{\mathbb{Z}}\, ^{ n_c}.$ In particular if
$\mathbb{X} $ is reduced to a single element, then the  index of
$\mathcal{T}$ is equal to the usual index of the semi-Fredholm
operator $\mathcal{T}.$ \eet

\bp From the properties   of the index \cite[Chapter 4]{CPY}, it
follows that two semi-Fredholm operators located in the same
connected component of the set of semi-Fredholm operators have the
same index. Moreover, as $ \mathcal{T} $ is continuous, the image
of a connected component of the topological space $ \mathbb{X},$
is included in a connected component of the set of semi-Fredholm
operators.  This shows that the index of a family of semi-Fredholm
operators  is well defined.  Moreover, it is clear that if $
\mathbb{X}$ is reduced to a single element, the  index of
$\mathcal{T}$ defined here is equal to the usual index of the
single semi-Fredholm operator $T.$

\ep

\bdf A continuous  family $\mathcal{K}$  from  $\mathbb{X}$ to
$L(X)$ is said to be compact if $\mathcal{K}_x$ is compact for all
$ x \in X.$

\edf

\bprop i) Let $\mathcal{T } \in  \mathcal{C}( \mathbb{X},$ {\bf
\scriptsize S}$\Phi(X))$ be a continuous  semi-Fredholm family and
let $ \mathcal{K}$ be a  continuous compact family.  Then $
\mathcal{T}+\mathcal{K} $ is a semi-Fredholm family and $
ind(\mathcal{T})= ind(\mathcal{T}+\mathcal{K}).$

ii) Let $\mathcal{S,T } \in  \mathcal{C}( \mathbb{X},$ {\bf
\scriptsize S}$\Phi(X))$ be two  semi-Fredholm families, such that
for all $ x \in \mathbb{X},  $  $\mathcal{S}_x$ and $\mathcal{T}_x
$ are both upper semi-Fredholm or both lower semi-Fredholm. Then
the family $\mathcal{ST}$ defined by $(\mathcal{ST})_x=
\mathcal{S}_x\mathcal{T}_x$ is a continuous semi-Fredholm family
and $ind(\mathcal{ST})= ind(\mathcal{S}) + ind(\mathcal{T}).$

\eprop

\bp This is clear from the usual properties of semi-Fredholm
operators.

\ep

\bdf Let $ \mathcal{S,T} $   be in $\in \mathcal{C}( \mathbb{X}, $
{\bf \scriptsize S}$\Phi(X)).$ We will say that $\mathcal{S}$ and
$\mathcal{T}$ are semi-Fredholm homotopic, if there exists a
continuous map $ h:  [0,1]\times X \rightarrow L(X)$ such that
$h(0,x)= \mathcal{S}_x, \, h(1,x)= \mathcal{T}_x $ and for all $
(t,x) \in [0,1]\times \mathbb{X},\, h(t,x)$ is a semi-Fredholm
operator. \edf

\bet Let $\mathcal{ S, T}$ be two semi-Fredholm homotopic elements
of \,$\mathcal{C}( \mathbb{X}, L(X)). $  Then $ind(\mathcal{T})=
ind(\mathcal{S}).$ \eet

\bp Since $\mathcal{S}$ and $\mathcal{T}$ are semi-Fredholm
homotopic, there exists a continuous map $ h: \mathbb{X} \times
[0, 1] \rightarrow Fred(X) $  such that such that
$h(0,x)=\mathcal{S}(x)$ and $h(1,x)=T(x)$ for all
$x\in\mathbb{X}.$ For a fixed $x \in \mathbb{X} ,$ the map $ h_x:
[0,1] \rightarrow $ {\bf \scriptsize S}$\Phi(X)),$ defined by
$h_x(t)= h(t,x)$ is a continuous path
 in  {\bf
\scriptsize S}$\Phi(X))$ linking  $\mathcal{S}_x$ to
$\mathcal{T}_x.$  Since $[0,1]$ is connected and since the
 usual index is constant on connected sets of  {\bf
\scriptsize S}$\Phi(X)$, then $ind(\mathcal{S}_x)=
 ind(\mathcal{T}_x).$  Therefore  $ind(\mathcal{S})= ind(\mathcal{T}).$
\ep

\bet  Let $\mathbb{X}$ be a compact topological space.  Then the
set $ \mathcal{C}( \mathbb{X}, \mathcal{SF}(X))$ of continuous
semi-Fredholm families is an open subset of the Banach algebra $
\mathcal{C}( \mathbb{X}, L(X))$ endowed with the uniform norm $
\norm{ \mathcal{T}}= sup_{x\in \mathbb{X}} \norm{\mathcal{T}_x}.$
\eet

\bp Let  $\mathcal{T} \in \mathcal{C}( \mathbb{X},$  {\bf
\scriptsize S}$\Phi(X))$),  then $\forall x \in \mathbb{X},
\exists \, \epsilon_x
> 0, $ such that $B(\mathcal{T}_x, \epsilon_x) \subset $ {\bf
\scriptsize S}$\Phi(X).$  We have $ \mathbb{X} \subset \bigcup
\limits _{x \in \mathbb{X}}{} \mathcal{T}^{-1}( B(\mathcal{T}_x,
\frac{\epsilon_{x_i}}{2})).$ Since $\mathbb{X}$ is compact, there
exists $ x_1,...x_n$ in $\mathbb{X}$ such that $\mathbb{X} \subset
\bigcup \limits_{i=1}^n T^{-1}(B(\mathcal{T}_{x_i},
\frac{\epsilon_{x_i}}{2})).$ Let $ \epsilon= min \{
\frac{\epsilon_{x_i}}{2} | 1\leq i \leq n \} $ the minimum of the
$ \frac{\epsilon_{x_i}}{2}, 1\leq i \leq n.$ Let $\mathcal{S} \in
\mathcal{C}( \mathbb{X}, L(X)) $ such that $ || \mathcal{T-S}|| <
\frac{\epsilon}{2}.$ If $ x \in \mathbb{X},$ then $ ||\mathcal{
T}_x-\mathcal{S}_x|| < \frac{\epsilon}{2}$ and there exists $ i,
1\leq i \leq n,$ such that $ x \in \mathcal{T}^{-1}(B(T_{x_i},
\epsilon_{x_i})).$ Then  $ || \mathcal{S}_x -\mathcal{T}_{x_i}||
\leq || \mathcal{S}_x-\mathcal{T}_x||  + ||\mathcal{T}_x-
\mathcal{T}_{x_i}|| < \epsilon/2 + \epsilon_{x_i}/2 \leq
\epsilon_{x_i}$ and so  $ \mathcal{S}_x$ is a semi-Fredholm
operator. Therefore $\mathcal{S} \in \mathcal{C}( \mathbb{X}, $
{\bf \scriptsize S}$\Phi(X))$ and $\mathcal{C}( \mathbb{X}, $ {\bf
\scriptsize S}$\Phi(X))$ is open in $ \mathcal{C}( \mathbb{X},
L(X)).$

\ep

\bet \label {index-continuity} Let $\mathbb{X}$ be a compact
topological  space. Then the index is a continuous locally
constant function  from $ \mathcal{C}( \mathbb{X}, $ {\bf
\scriptsize S}$\Phi(X))$into the group $\mathbb{Z}^{
n_\textit{\textbf{c}}}.$

\eet

\bp As the set   {\bf \scriptsize S}$\Phi(X))$ is an open subset
of $L(X),$  if $\mathcal{T} \in \mathcal{C}( \mathbb{X},$  {\bf
\scriptsize S}$\Phi(X))$), then $\forall x \in \mathbb{X}, \exists
\, \epsilon_x
> 0, $ such that $B(\mathcal{T}_x, \epsilon_x) \subset $ {\bf
\scriptsize S}$\Phi(X)$. Then the index is constant on
$B(\mathcal{T}_x, \epsilon_x),$ because $B(\mathcal{T}_x,
\epsilon_x)$ is connected. We have $ \mathbb{X} \subset \bigcup
\limits _{x \in \mathbb{X}}{}\mathcal{T}^{-1}( B(\mathcal{T}_x,
\frac{\epsilon_{x}}{2})).$ Since $\mathbb{X}$ is compact, there
exists $ x_1,...x_n$ in $\mathbb{X}$ such that $\mathbb{X} \subset
\bigcup \limits_{i=1}^n\mathcal{ T}^{-1}(B(T_{x_i},
\frac{\epsilon_{x_i}}{2})).$ Let $ \epsilon= min \{
\frac{\epsilon_{x_i}}{2} | 1\leq i \leq n \} $ the minimum of the
$ \frac{\epsilon_{x_i}}{2}, 1\leq i \leq n,$ and let $ S \in
\mathcal{C}( \mathbb{X},$  {\bf \scriptsize S}$\Phi(X))$ such that
$ || \mathcal{T-S}|| < \frac{\epsilon}{2}.$ If $ x \in
\mathbb{X},$ then $ || \mathcal{T}_x-\mathcal{S}_x|| <
\frac{\epsilon}{2}$ and there exists $ i, 1\leq i \leq n,$ such
that $ x \in\mathcal{ T}^{-1}(B(T_{x_i},
\frac{\epsilon_{x_i}}{2})).$ Then $ || \mathcal{S}_x
-\mathcal{T}_{x_i}|| \leq || \mathcal{S}_x-\mathcal{T}_x|| +
||\mathcal{T}_x- \mathcal{T}_{x_i}|| < \epsilon/2 +
\epsilon_{x_i}/2 \leq \epsilon_{x_i}.$ So  $ind(\mathcal{S}_x)=
ind(\mathcal{T}_{x_i})= ind(\mathcal{T}_x).$ Hence the index is a
locally constant function, in particular it is a continuous
function. \ep

\section {upper semiregularities of continuous families of bounded linear
operators}

In this part we consider various regularities of continuous
families of  bounded linear operators, parametrized by compact
topological space $ \mathbb{X}.$ These regularities had been
considered in \cite{MM} and \cite{P8}, in the case when  $
\mathbb{X}$ is reduced to a single element.

 \bdf Let $ \mathbb{X}$ be a compact topological space. We define
 the sets  $ \mathcal{R}_{i}, 1\leq i \leq 16,$ as follows:

\par\noindent   $ {\bf{\mathcal{R}}_{1}} =\{ \mathcal{T} \in \mathcal{C}( \mathbb{X},
L(X)) \mid \forall x \in \mathbb{X}, \mathcal{T}_x$ is onto\}  \\
\noindent $ {\bf{\mathcal{R}}_{2}} = \{ \mathcal{T} \in
\mathcal{C}( \mathbb{X}, L(X)) \mid \forall x \in \mathbb{X},
\mathcal{T}_x \in \Phi_{-} (X)$ and
$ \delta(\mathcal{T}_x) < \infty\}. $\\
  \noindent ${\bf{\mathcal{R}}_{3}} = \{\mathcal{T} \in \mathcal{C}( \mathbb{X},
L(X)) \mid \forall x \in \mathbb{X}, :  \delta(\mathcal{T}_x) <
\infty $ and $R(\mathcal{T}{_x}^{\delta(\mathcal{T}_x)})$
   is closed\}.\\
\noindent ${\bf{\mathcal{R}}_{4}} =  \{ \mathcal{T} \in
\mathcal{C}( \mathbb{X}, L(X)) \mid \forall x \in \mathbb{X},
\mathcal{T}_x \in    \Phi_{-} (X)\}. $\\
\noindent ${\bf{\mathcal{R}}_{5}} = \{ \mathcal{T} \in
\mathcal{C}( \mathbb{X}, L(X)) \mid \forall x \in \mathbb{X},
\delta_{e}(\mathcal{T}_x) < \infty$ and
  $R(\mathcal{T}^{\delta_{e}(\mathcal{T}_x)})$ is closed \}.\\
\noindent ${\bf{\mathcal{R}}_{6}} = \{ T \in \mathcal{C}(
\mathbb{X}, L(X))
\mid \forall x \in \mathbb{X}, \mathcal{T}_x$  is bounded below \}.\\
\noindent ${\bf{\mathcal{R}}_{7}} =\{\mathcal{T} \in \mathcal{C}(
\mathbb{X}, L(X)) \mid \forall x \in \mathbb{X}, \mathcal{T}_x \in
\Phi_{+} (X)$
 and   $    a(\mathcal{T}_x) < \infty\}. $\\
\noindent  ${\bf{\mathcal{R}}_{8}} = \{  \mathcal{T} \in
\mathcal{C}( \mathbb{X}, L(X)) \mid \forall x \in \mathbb{X},
a(\mathcal{T}_x) < \infty $ and
   $R(\mathcal{T}{_x}^{a(\mathcal{T}_x)+1}) $ \,\, is closed \}. \\
\noindent ${\bf{\mathcal{R}}_{9}} =  \{ \mathcal{T} \in
\mathcal{C}( \mathbb{X}, L(X))
\mid \forall x \in \mathbb{X}, \mathcal{T}_x \in    \Phi_{+} (X)\}. $   \\
\noindent  ${\bf{\mathcal{R}}_{10}} = \{ \mathcal{T} \in
\mathcal{C}( \mathbb{X}, L(X)) \mid \forall x \in \mathbb{X},
a_{e}(\mathcal{T}_x) < \infty $
 and $R(\mathcal{T}^{a_{e}(\mathcal{T}_x)+1})$ is closed \}.
\\
\noindent ${\bf{\mathcal{R}}_{11}} = \{ \mathcal{T} \in
\mathcal{C}( \mathbb{X}, L(X)) \mid \forall x \in \mathbb{X},
N(\mathcal{T}_x)\subset R^{\infty}(\mathcal{T}_x)$ and
$R(\mathcal{T}_x)$ is closed \},where  $R^{\infty}(\mathcal{T}_x)=
\bigcap \limits _{n \in
\mathbb{N}} R(\mathcal{T}_x^n).$  \\
\noindent ${\bf{\mathcal{R}}_{12}} = \{ \mathcal{T} \in
\mathcal{C}( \mathbb{X}, L(X)) \mid \forall x \in \mathbb{X},
N(\mathcal{T}_x){\subset}_{e} R^{\infty}(\mathcal{T}_x)$ and
$R(\mathcal{T}_x)$ is
closed \}.\\
\noindent ${\bf{\mathcal{R}}_{13}} = \{ \mathcal{T} \in
\mathcal{C}( \mathbb{X}, L(X)) \mid \forall x \in \mathbb{X},
\exists \,\,  d_x \in  \mathbb{N}:
R(\mathcal{T}_x)+N(\mathcal{T}_x^{d_x})=
R(\mathcal{T}_x)+N^{\infty}(\mathcal{T}_x) $ and $
R(\mathcal{T}_x^{d_x+1})$
 is closed \},  where  $N^{\infty}(\mathcal{T}_x)= \bigcup \limits _{n \in
\mathbb{N}} N(\mathcal{T}_x^n).$\\
\noindent ${\bf{\mathcal{R}}_{14}} = \{ \mathcal{T} \in
\mathcal{C}( \mathbb{X}, L(X)) \mid \forall x \in \mathbb{X},
\forall n\in \mathbb{ N}, k_{n}(\mathcal{T}_x)
< \infty $  and $R(\mathcal{T}_x)$ is closed\}.\\
\noindent ${\bf{\mathcal{R}}_{15}}= \{ \mathcal{T} \in
\mathcal{C}( \mathbb{X}, L(X)) \mid \forall x \in \mathbb{X},
\exists \,\,d_x \in \mathbb{ N}\,\,: k_{n}(\mathcal{T}_x)<\infty $ if $ n\geq d_x $ and $R(\mathcal{T}_x^{d_x+1})$ is closed \}. \\
\noindent ${\bf{\mathcal{R}}_{16}}= \{\mathcal{T} \in \mathcal{C}(
\mathbb{X}, L(X)) \mid \forall x \in \mathbb{X}, \,\exists \,
d_x\in \mathbb{ N}: \mathcal{T}_x  $ is of  topological uniform
descent for  $ n \geq d_x \}.$

\edf

 \noindent We have ${\bf{\mathcal{R}}_{1}} \subset {\bf{\mathcal{R}}_{2}}\subset {\bf{\mathcal{R}}_{3}}
     \subset  {\bf{\mathcal{R}}_{3}}\cup {\bf{\mathcal{R}}_{4}} \subset {\bf{\mathcal{R}}_{5}}
     \subset {\bf{\mathcal{R}}_{13}}$,
     ${\bf{\mathcal{R}}_{6}}\subset {\bf{\mathcal{R}}_{7}} \subset  {\bf{\mathcal{R}}_{8}} \subset
     {\bf{\mathcal{R}}_{8}}\cup  {\bf{\mathcal{R}}_{9}} \subset  {\bf{\mathcal{R}}_{10}}\subset {\bf{\mathcal{R}}_{13}},$
 ${\bf{\mathcal{R}}_{11}} \subset {\bf{\mathcal{R}}_{12}} \subset {\bf{\mathcal{R}}_{13}}
 \subset{\bf{\mathcal{R}}_{13}}\cup {\bf{\mathcal{R}}_{14}}\subset
 {\bf{\mathcal{R}}_{15}} \subset
 {\bf{\mathcal{R}}_{16}}$. The elements  of $ {\bf{\mathcal{R}}_{11}}$ and  $ {\bf{\mathcal{R}}_{12}}$
 are called semi-regular and essentially semi-regular continuous families.
 The elements of $ {\bf{\mathcal{R}}_{13}}$ are called  quasi-Fredholm
continuous families  and  the elements of  $
{\bf{\mathcal{R}}_{16}}$ are called  continuous families of
topological uniform descent.

\noindent The inclusion ${\bf{\mathcal{R}}_{5}}\cup
{\bf{\mathcal{R}}_{10}}\subset
 {\bf{\mathcal{R}}_{13}}$ is a consequence of  \cite[Proposition 3.4]{P8}. \\

\bet The sets $ {\bf {\mathcal{R}_{i}}},$ for  $1\leq i \leq 16, $
are regularities in the Banach algebra  $\mathcal{C}( \mathbb{X},
L(X)). $

\eet

\bp From \cite{MM}  and \cite[Theorem 4.3]{P8}, it follows that
for each single element $x$ of $\mathbb{X}$, the set $
{\bf{\mathcal{R}}_{i,x}}= \{ \mathcal{T}_x \mid \mathcal{T} \in
{\bf{\mathcal{R}}_{i}} \},$  $1\leq i \leq 16 $ is a regularity.
Thus  ${\bf{\mathcal{R}}_{i}},$  $1\leq i \leq 16, $ are
regularities.

\ep

\bcor \label{smt}   Let  $\mathcal{\mathcal{T}} \in \mathcal{C}(
\mathbb{X}, L(X))$  and $f$ an analytic function   in a
 neighborhood  of the usual spectrum $\sigma(\mathcal{T})$  of   $\mathcal{T}$ which
 is non-constant on any connected component of  $\sigma(\mathcal{T})$.
Then
$f(\sigma_{\bf{\mathcal{R}}_{i}}(\mathcal{T}))=\sigma_{\bf{\mathcal{R}}_{i}}(f(\mathcal{T})),
1\leq i \leq 16. $ \ecor

\bp  From \cite{MM}  and \cite[Theorem 4.3]{P8}, it follows that
for each single element $x$ of $\mathbb{X}$, the set $
{\bf{\mathcal{R}}_{i,x}}= \{ \mathcal{T}_x \mid \mathcal{T} \in
{\bf{\mathcal{R}}_{i}} \}$ for $1\leq i \leq 16, $ is a
regularity. As  $ \sigma_(\mathcal{\mathcal{T}})= \bigcup \limits
_{x \in \mathbb{X}}{} \sigma (\mathcal{\mathcal{T}}_x), 1\leq i
\leq 16, $ then for all $ x \in \mathbb{X},f $ is an analytic
function in a
 neighborhood  of the usual spectrum $\sigma(\mathcal{T}_x)$ of
 $x.$

From \cite{MM}, it follows that
$f(\sigma_{{\bf{\mathcal{R}}_{i}}}(\mathcal{T}_x))=\sigma_{{\bf{\mathcal{R}}_{i}}}(f(\mathcal{T}_x)),
1\leq i \leq 15$      and from \cite[Theorem 4.3]{P8}, it follows
that
$f(\sigma_{{\bf{\mathcal{R}}_{16}}}(\mathcal{T}_x))=\sigma_{{\bf{\mathcal{R}}_{16}}}
(f(\mathcal{T}_x)),$  for all $ x \in \mathbb{X}.$ As  $
\sigma_{{\bf{\mathcal{R}}_{i}}}(\mathcal{\mathcal{T}})= \bigcup
\limits _{x \in \mathbb{X}}{}
\sigma_{{\bf{\mathcal{R}}_{i,x}}}(\mathcal{\mathcal{T}}_x), $ then
$f(\sigma_{{\bf{\mathcal{R}}_{i}}}(\mathcal{\mathcal{T}}))=\sigma_{{\bf{\mathcal{R}}_{i}}}(f(\mathcal{\mathcal{T}})),
1\leq i \leq 16.  $ \ep

\bdf A continuous  family $\mathcal{F}$  from  $\mathbb{X}$ to
$L(X)$ is said to be of finite rank if $\mathcal{F}_x$ is of
finite rank  for all $ x \in\mathbb{ X}.$

\edf

\bcor
 Let  $\mathcal{T} \in {\bf{\mathcal{R}}_{i}},$ for $ i \in \{
 4,5,9,10,12,13,14,15 \}.$  If  $\mathcal{F}$  is a continuous family  of finite
 rank, then $ \mathcal{T+ F} \in {\bf{\mathcal{R}}_{i}}.$
\ecor

\bp From \cite{MM}, it follows that for each single element $x$ of
$\mathbb{X}$, the regularity  $ {\bf{\mathcal{R}}_{i,x}}, i \in \{
4,5,9,10,11,12,13,14,15 \},$ is stable under the  finite rank
perturbation $\mathcal{F}_x$ . Thus $ \mathcal{T+ F} \in
{\bf{\mathcal{R}}_{i}}.$ \ep

\bremark  A continuous family $ \mathcal{T}$ in $\mathcal{C}(
\mathbb{X}, L(X))$  is said to be a Fredholm family if for all $ x
\in \mathbb{X}, \mathcal{T}_x$ is a  Fredholm operator. We observe
that the set $\Phi(X)$ of continuous Fredholm families  is equal
to ${\bf{\mathcal{R}}_{4}} \cap {\bf{\mathcal{R}}_{9}},$ where
${\bf{\mathcal{R}}_{4}}$  is the regularity of continuous families
of upper semi-Fredholm operators and  ${\bf{\mathcal{R}}_{9}}$ is
the regularity of continuous families of lower semi-Fredholm
operators. Thus $\Phi(X)$ is a regularity and
$\sigma_{\Phi(X)}(\mathcal{T})=
\sigma_{\bf{\mathcal{R}}_{4}}(\mathcal{T}) \cup
\sigma_{\bf{\mathcal{R}}_{9}}(\mathcal{T}).$

\eremark

Now, we extend some known results \cite[Theorem 4.3.2, Theorem
4.3.3, Corollary 1.3.4, Corollary 1.3.5, Corollary 1.3.6] {CPY} in
the case of a single semi-Fredholm operators to the case of
continuous families of semi-Fredholm operators. As their proofs
are  straightforwardly  obtained from the case of a single
operator, we give them without proof.

\bet

\begin{enumerate}

\item  Let  $\mathcal{\mathcal{T}} \in \mathcal{C}( \mathbb{X},
L(X)).$ Then $\mathcal{\mathcal{T}}$ is in
${\bf{\mathcal{R}}_{4}}$  if and only if there exists a family of
bounded linear operators $ (U_x)_{x \in \mathbb{X}}, $  a family
of bounded linear compact operators $ (K_x)_{x \in \mathbb{X}}, $
a family of bounded projections   $ (P_x)_{x \in \mathbb{X}}, $
such that for all $ x \in \mathbb{X},$  $ P_x$  is defined from
$X$ onto $N( \mathcal{T}_x)$ and $ \mathcal{T}_xU_x= I + K_x.$

\item  Let  $\mathcal{\mathcal{T}} \in \mathcal{C}( \mathbb{X},
L(X)).$ Then $\mathcal{\mathcal{T}}$ is in
${\bf{\mathcal{R}}_{9}}$ if and only if there exists a family of
bounded linear operators $ (V_x)_{x \in \mathbb{X}}, $  a family
of bounded linear compact operators $ (K_x)_{x \in \mathbb{X}}, $
a family of bounded projections   $ (P_x)_{x \in \mathbb{X}}, $
such that for all $ x \in \mathbb{X},$  $ P_x$  is defined from
$X$ onto $R( \mathcal{T}_x)$ and $ V_x \mathcal{T}_x= I + K_x.$

\item  Let  $\mathcal{\mathcal{T}}_1,  \mathcal{\mathcal{T}}_2$ in
$\mathcal{C}( \mathbb{X}, L(X))$ such that
$\mathcal{\mathcal{T}}_1\mathcal{\mathcal{T}}_2  \in
{\bf{\mathcal{R}}_{9}, }$   then  $\mathcal{\mathcal{T}}_2  \in
{\bf{\mathcal{R}}_{9}.}$

\item  Let  $\mathcal{\mathcal{T}}_1,  \mathcal{\mathcal{T}}_2$ in
$\mathcal{C}( \mathbb{X}, L(X))$ such that
$\mathcal{\mathcal{T}}_1\mathcal{\mathcal{T}}_2  \in
{\bf{\mathcal{R}}_{4}, }$   then  $\mathcal{\mathcal{T}}_1 \in
{\bf{\mathcal{R}}_{4}.}$

\item  Let  $\mathcal{\mathcal{T}}_1,  \mathcal{\mathcal{T}}_2$ in
$\mathcal{C}( \mathbb{X}, L(X))$ such that
$\mathcal{\mathcal{T}}_1\mathcal{\mathcal{T}}_2  \in \Phi(X),$
then $\mathcal{\mathcal{T}}_1  \in {\bf{\mathcal{R}}_{4}}$ and
$\mathcal{\mathcal{T}}_2 \in {\bf{\mathcal{R}}_{9}.}$

\end{enumerate}

\eet

\section { Upper and lower semiregularities  of continuous families of bounded linear
operators}

Following \cite{MU},  the axioms i) and ii) of regularities in
Definition \ref{REG},  are splitted in two  halves, each of them
implying a one-way spectral mapping theorem and giving raises to
the notions of upper and lower semiregularities. We first consider
lower semiregularities.

In all what follows,  $ \mathbb{X} $ will be a compact topological
space. If
 ${\bf {\mathcal{R}}}$
  is an nonempty subset  of \, $\mathcal{C}( \mathbb{X},L(X)),$  then for all $x \in \mathbb{X},$ we set
   ${\bf {\mathcal{R}}}_x= \{\mathcal{T}_x \mid \mathcal{T} \in {\bf {\mathcal{R}}}
   \}, x \in \mathbb{X}.$ \\

{ \bf \underline{A- Upper semiregularities} }\\

\bdf
 Let $A$ be a unital Banach algebra, with unit $\mathbf{e}.$
 A non-empty subset  ${\bf{\mathcal{R}}}$ of $A$ is called a lower
semiregularity  if it satisfies the following conditions:\\
 \noindent i) If $ a \in A$ and  $a^n \in {\bf R},$  $n\geq 1,$ then  $
a \in {\bf{\mathcal{R}}} .$\\
 \noindent ii) If $ a,b,c,d \in A $  are mutually commuting
 elements of $A$
satisfying  $ ac+bd= \mathbf{e}, $ such that  $ab \in {\bf
{\mathcal{R}}}, $ then
 $a,b \in {\bf {\mathcal{R}}} $ .
\edf

\noindent If  ${\bf {\mathcal{R}}}$  is a lower semiregularity of
$A$ and $ a \in A,$ we define the  $\mathbf{R}$-spectrum of $a$ by
$ \sigma_{\mathbf{R}}(a) = \{ \lambda \in {\bf C}/ a-\lambda
\mathbf{e} \notin \mathbf{R}\} .$\\

Recall the following result from \cite[Lemma 2]{MU}

\begin{lem}
Let $A$ be a unital Banach algebra with unit $\mathbf{e}$ and let
 ${\bf {\mathcal{R}}}$  a lower semiregularity of $A.$ then:

 \begin{enumerate}

\item  $e \in  {\bf {\mathcal{R}}}.$ \item  $Inv(A) \subset {\bf
{\mathcal{R}}},$ where $Inv(A)$ is the set of invertible elements
of $A.$ \item  If $a \in {\bf {\mathcal{R}}}, b \in Inv (A)$ and
$ab = ba$, then $ab \in {\bf {\mathcal{R}}}.$
 \item  $\sigma_{{\bf {\mathcal{R}}}}(a) \subset  \sigma(a).$
\item (Translation property) $\sigma_{{\bf {\mathcal{R}}}}(a +
\lambda \mathbf{e})= \lambda \mathbf{e} + \sigma_{{\bf
{\mathcal{R}}}}(a),$ for all $\lambda \in \mathbb{C}.$

\end{enumerate}

\end{lem}

Let $ \mathbb{X} $ be a compact topological space, let
 ${\bf {\mathcal{R}}}$
  be an empty subset  of \, $\mathcal{C}( \mathbb{X},L(X)).$ For all $x \in \mathbb{X},$ let
   ${\bf {\mathcal{R}}}_x= \{\mathcal{T}_x \mid \mathcal{T} \in {\bf {\mathcal{R}}}
   \}, x \in \mathbb{X}.$ We assume that   ${\bf {\mathcal{R}}}$  satisfies the following two
conditions, summarized as condions (*) :

i) $ \forall x \in \mathbb{X}, \forall \mathcal{T} \in
\mathcal{C}( \mathbb{X},L(X)),  \mathcal{T}_x  \in {\bf
{\mathcal{R}}}_x \Rightarrow \mathcal{T} \in {\bf {\mathcal{R}}}.$

ii)  $ \forall x \in \mathbb{X},  {\bf {\mathcal{R}}}_x $ is a
lower semiregularity  of $L(X).$

Using the properties of  lower semiregularities, we  obtain a one
way spectral mapping theorem.

\bet \label{LR}
 Let $ \mathbb{X} $ be a compact topological space and let
 ${\bf {\mathcal{R}}}$ be a nonempty subsetof \, $\mathcal{C}( \mathbb{X},L(X)),$
satisfying the conditions (*).   Then:

\begin{enumerate}

  \item  ${\bf {\mathcal{R}}}$ is a lower  semiregularity.

  \item If $f$ an analytic function   in a
 neighborhood  of the usual spectrum $\sigma(\mathcal{T})$  of an element   $\mathcal{T}  \in \mathcal{C}( \mathbb{X},L(X)),$ which
 is non-constant on any connected component of
 $\sigma(\mathcal{T})$,
then $f(\sigma_{\bf{\mathcal{R}}}(\mathcal{T}))
\subset\sigma_{\bf{\mathcal{R}}}(f(\mathcal{T})). $

\end{enumerate}

\eet

\bp  1. Assume  for all  $ x \in \mathbb{X}, $ the set ${\bf
{\mathcal{R}}}_x$ is a lower
  semiregularity of $ L(X).$ Let us show that ${\bf{\mathcal{R}}}$
  is a regularity.

For if  $\mathcal{T} \in  {\bf{\mathcal{R}}}$  and $ n\geq 1$ such
that  $\mathcal{T}^n \in {\bf{\mathcal{R}}},$ then for all $x \in
\mathbb{X},$ we have $(\mathcal{T}^n)_x = (\mathcal{T}_x)^n \in
{\bf{\mathcal{R}_x}}.$ As ${\bf{\mathcal{R}_x}}$  is a lower
semiregularity, then for all $ x \in \mathbb{X},$ we have
$\mathcal{T}_x \in {\bf{\mathcal{R}_x}}.$ Since
${\bf{\mathcal{R}}}$ satisfies the conditions (*), then
$\mathcal{T}  \in  {\bf{\mathcal{R}}}.$

Now if  $\mathcal{S,T, U,V,} $ are commuting elements of
${\bf{\mathcal{R}}}$ satisfying $\mathcal{ US + VT = I},$ such
that $ \mathcal{ST} \in  {\bf{\mathcal{R}}}.$  Then for all $ x
\in \mathbb{X},$ we have  $\mathcal{U}_x \mathcal{S}_x +
\mathcal{V}_x\mathcal{T}_x =I$ and $\mathcal{S}_x,\mathcal{T}_x,
\mathcal{U}_x,\mathcal{V}_x, $ are commuting elements of
${\bf{\mathcal{R}_x}}.$ As we have   for all $x \in \mathbb{X},$ $
\mathcal{S}_x\mathcal{T}_x = (\mathcal{ST})_x \in
{\bf{\mathcal{R}_x}},$  and since ${\bf{\mathcal{R}}}_x$ is a
lower semiregularity, then $\mathcal{S}_x $ and  $\mathcal{T}_x
\in {\bf{\mathcal{R}}}_x.$ As $ {\bf{\mathcal{R}}}$ satifies
conditions (*), then $\mathcal{S} $ and $\mathcal{T}$  are in
${\bf{\mathcal{R}}}.$
Consequently ${\bf {\mathcal{R}}}$ is a lower  semiregularity in  $\mathcal{C}( \mathbb{X},L(X)).$ \\

\noindent 2.  Let us show that  $
\sigma_{\bf{\mathcal{R}}}(\mathcal{T})= \bigcup \limits _{x \in
\mathbb{X}}{}\sigma_{\bf{\mathcal{R}}_x}( \mathcal{T}_x).$ If
$\lambda \notin \sigma_{\bf{{\bf {\mathcal{R}}}}}(\mathcal{T}),$
then $ \mathcal{T} - \lambda \mathcal{I} \in {\bf{{\bf
{\mathcal{R}}}}}.$ Therefore  for all $ x \in \mathbb{X}$, we have
$ \mathcal{T}_x - \lambda I \in {\bf{{\bf {\mathcal{R}}}}}_x.$
Hence $\lambda \notin \sigma_{\bf{\mathcal{R}}_x}(
\mathcal{T}_x).$ Conversely if  for all  $ x \in \mathbb{X},
\lambda \notin \sigma_{\bf{\mathcal{R}}_x}( \mathcal{T}_x), $ then
for all  $ x \in \mathbb{X}, \mathcal{T}_x - \lambda I \in
{\bf{\mathcal{R}}_x}.$ As $ {\bf{\mathcal{R}}}$ satisfies the
conditions (*), then $ \mathcal{T} - \lambda \mathcal{I} \in
{\bf{{\bf {\mathcal{R}}}}}.$ So $ \sigma_{\bf{{\bf
{\mathcal{R}}}}}(\mathcal{T})= \bigcup \limits _{x \in
\mathbb{X}}{}\sigma_{\bf{\mathcal{R}}_x}( \mathcal{T}_x).$

As  $ \sigma(\mathcal{T})=   \bigcup \limits _{x \in
\mathbb{X}}{}\sigma( \mathcal{T}_x), $ then  for all $ x \in
\mathbb{X},$ $f$ is an analytic function   in a
 neighborhood  of the usual spectrum $\sigma(\mathcal{T}_x)$  of
 $\mathcal{T}_x.$ From \cite[Theorem 4]{MU},  we have
$f(\sigma_{\bf{{\bf {\mathcal{R}}}}_x}(\mathcal{T}_x))
\subset\sigma_{\bf{{\bf {\mathcal{R}}}_x}}(f(\mathcal{T}_x)). $
Then $f(\sigma_{\bf{{\bf {\mathcal{R}}}}}(\mathcal{T}))= f(
\bigcup \limits_{x \in \mathbb{X}}{}\sigma_{\bf{\mathcal{R}}_x}(
\mathcal{T}_x))= \bigcup \limits _{x \in
\mathbb{X}}{}f(\sigma_{\bf{\mathcal{R}}_x}( \mathcal{T}_x))
\subset \bigcup \limits _{x \in
\mathbb{X}}{}\sigma_{\bf{\mathcal{R}}_x}(f( \mathcal{T}_x))
=\sigma_{\bf{{\bf {\mathcal{R}}}}}(f(\mathcal{T}),$ because by the
properties of the functional calculus,  we have $f(
\mathcal{T}_x)= f( \mathcal{T})_x.$

 \ep

\bex  Let  $ (m_x)_{x \in \mathbb{X}}$ be a sequence in $
\mathbb{N}^{\mathbb{X}}.$ From \cite{MU} and Proposition \ref{LR},
we deduce that the following sets are lower semiregularities of
$\mathcal{C}( \mathbb{X}, L(X)).$
\begin{enumerate}

\item $  \mathcal{LSR}_1=  \{\mathcal{T} \in \mathcal{C}(
\mathbb{X}, L(X)) : \forall x \in \mathbb{X}, \,\,
\alpha(\mathcal{T}_x) \leq m_x  \}$

\item $ \mathcal{LSR}_2= \{\mathcal{T} \in \mathcal{C}(
\mathbb{X}, L(X)) : \forall x \in \mathbb{X},\,
\beta(\mathcal{T}_x) \leq  m_x \}$

\item $\mathcal{LSR}_3=  \{ \mathcal{T} \in \mathcal{C}(
\mathbb{X}, L(X)) :\forall x \in \mathbb{X}, \,
\alpha(\mathcal{T}_x) \leq m_x$ and $R(\mathcal{T}_x)$ is closed
\}

\end{enumerate}

\eex

{ \bf \underline{B- Upper semiregularities} }\\

 We assume that   ${\bf {\mathcal{R}}}$  satisfies the following two
conditions, summarized as conditions (**) :

i) $ \forall x \in \mathbb{X}, \forall \mathcal{T} \in
\mathcal{C}( \mathbb{X},L(X)),  \mathcal{T}_x  \in {\bf
{\mathcal{R}}}_x \Rightarrow \mathcal{T} \in {\bf {\mathcal{R}}}.$

ii)  $ \forall x \in \mathbb{X},  {\bf {\mathcal{R}}}_x $ is an
upper semiregularity  of $L(X).$

\bdf \label{UR}
 Let $A$ be a unital Banach algebra, with unit $\mathbf{e}.$
  A non-empty subset ${\bf {\mathcal{R}}}$ of $A$ is called an  upper
semiregularity  if it satisfies the following conditions:\\
 \noindent i) If $ a \in A$ and  $a \in {\bf {\mathcal{R}}},$  then for all $n\geq 1,$ then  $
a^n \in {\bf {\mathcal{R}}}.$\\
 \noindent ii) If $ a,b,c,d \in L(X) $  are mutually commuting
 elements of $A$
satisfying  $ ac+bd= \mathbf{e}, $ then  if $a, b \in {\bf
{\mathcal{R}}}, $ then
 $ab \in {\bf {\mathcal{R}}}$ .\\
\noindent  iii) $ {\bf {\mathcal{R}}} $ contains a neighbourhood
of the unit element $\mathbf{e}.$

\edf

\noindent If  ${\bf {\mathcal{R}}}$  is a lower semiregularity of
$A$ and $ a \in A,$ we define the  ${\bf {\mathcal{R}}}$-spectrum
of $a$ by $ \sigma_{{\bf {\mathcal{R}}}}(a) = \{ \lambda \in {\bf
C}/ a-\lambda \mathbf{e}\notin
{\bf {\mathcal{R}}}\} .$  \\

Using the properties of  upper semiregularities, we  obtain a one
way spectral mapping theorem.  Here $ Inv(\mathcal{C}(
  \mathbb{X},L(X))$ and $Inv(L(X)),$ denote respectively the set of invertible
  elements of $\mathcal{C}(\mathbb{X},L(X)$ and $L(X)$ .

\bet \label{UR1}
 Let $ \mathbb{X} $ be a compact topological space and let
 $\mathcal{R}$
  be a nonempty subset of \, $\mathcal{C}( \mathbb{X},L(X)),$ containing an open
   neighborhood of the identity $\mathcal{I}$ and satisfying the conditions (**). Then:

\begin{enumerate}

  \item$\mathcal{R}$ is an upper
  semiregularity of \, $\mathcal{C}( \mathbb{X},L(X)).$
  \item Assume that: $$\forall  x \in \mathbb{X},\,  { \text if}\,
\mathcal{T}_x \in {\bf\mathcal{R}}_x \cap Inv(L(X)),
  then \,\, (\mathcal{T}_x)^{-1} \in  {\bf{\mathcal{R}}}. \, (***)$$
Then if $f$ an analytic function   in a
 neighborhood  of the usual spectrum $\sigma(\mathcal{T})$  of   $\mathcal{T}$ which
 is non-constant on any connected component of
 $\sigma(\mathcal{T})$,
then $ \sigma_{\bf{\mathcal{R}}}(f(\mathcal{T}) \subset
f(\sigma_{\bf{\mathcal{R}}}(\mathcal{T})). $

\end{enumerate}

\eet

\bp  \begin{enumerate}  \item  Assume that  for all  $ x \in
\mathbb{X}, $ the set ${\bf {\mathcal{R}}}_x $ is an upper
  semiregularity of $ L(X).$    Let us show that ${\bf{\mathcal{R}}}$
  is an upper  regularity.

For if  $\mathcal{T} \in  {\bf{\mathcal{R}}},$  then for all $x
\in \mathbb{X},$ we have $\mathcal{T}_x  \in
{\bf{\mathcal{R}_x}}.$ As ${\bf{\mathcal{R}_x}}$  is an upper
semiregularity, then for all $ x \in \mathbb{X},$ we have
$(\mathcal{T}_x)^n=   (\mathcal{T}^n)_x \in {\bf{\mathcal{R}_x}}.$
Since ${\bf{\mathcal{R}}}$ satisfies the conditions (**), then
$\mathcal{T}^n  \in  {\bf{\mathcal{R}}}.$

 Now if  $\mathcal{S,T, U,V,} $ are commuting elements of
${\bf{\mathcal{R}}}$ satisfying $\mathcal{ US + VT = I},$ then for
all $ x \in \mathbb{X},$ we have  $\mathcal{U}_x \mathcal{S}_x +
\mathcal{V}_x\mathcal{T}_x =I$ and $\mathcal{S}_x,\mathcal{T}_x,
\mathcal{U}_x,\mathcal{V}_x, $ are commuting elements of
${\bf{\mathcal{R}_x}}.$  As ${\bf{\mathcal{R}_x}}$  is an upper
semiregularity, then for all $x \in \mathbb{X},$ we have $
\mathcal{S}_x\mathcal{T}_x = (\mathcal{ST})_x \in
{\bf{\mathcal{R}_x}}.$ Since ${\bf{\mathcal{R}}}$ satisfies the
conditions (**), then $\mathcal{ST}  \in  {\bf{\mathcal{R}}}.$\\
As we know by hypothesis that ${\bf {\mathcal{R}}}$  contains an
open
   neighborhood of the identity $\mathcal{I},$ then
 ${\bf {\mathcal{R}}}$ is an upper semiregularity in $\mathcal{C}(
\mathbb{X},L(X)).$

\item  As for all $x \in \mathbb{X}, {\bf{\mathcal{R}_x}}$  is an
upper semiregularity  satisfying (***), then from \cite[Theorem
20]{MU},  we have $
\sigma_{\bf{\mathcal{R}_x}}(f(\mathcal{T}_x))\subset
f(\sigma_{\bf{\mathcal{R}}_x}(\mathcal{T}_x)).$  Then $
\sigma_{\bf{\mathcal{R}}}(f(\mathcal{T}))=   \bigcup \limits _{x
\in \mathbb{X}}{}\sigma_{\bf{\mathcal{R}}_x}( f(\mathcal{T})_x)=
\bigcup \limits _{x \in \mathbb{X}}{}\sigma_{\bf{\mathcal{R}}_x}(
f(\mathcal{T}_x)),$
  because by the properties of the functional calculus,
we have   $f( \mathcal{T})_x = f( \mathcal{T}_x).$ So $
\sigma_{\bf{\mathcal{R}}}(f(\mathcal{T})) \subset \bigcup \limits
_{x \in \mathbb{X}}{}f(\sigma_{\bf{\mathcal{R}}_x}(\mathcal{T}_x))
$ and $ \sigma_{\bf{\mathcal{R}}}(f(\mathcal{T})) \subset
f(\sigma_{\bf{\mathcal{R}}}(\mathcal{T}))),$ because $\bigcup
\limits _{x \in
\mathbb{X}}{}\sigma_{\bf{\mathcal{R}}_x}(\mathcal{T}_x)=
\sigma_{\bf{\mathcal{R}}}(\mathcal{T})$ (see the proof of Theorem
\ref{LR}) and $ f(\bigcup \limits _{x \in
\mathbb{X}}{}\sigma_{\bf{\mathcal{R}}_x}(\mathcal{T}_x)) = \bigcup
\limits _{x \in
\mathbb{X}}{}f(\sigma_{\bf{\mathcal{R}}_x}(\mathcal{T}_x)). $

\end{enumerate}

 \ep

\bex  Let  $ (m_x)_{x \in \mathbb{X}}$ be a sequence in $
\mathbb{N}^{\mathbb{X}}.$ From  Theorem \ref{UR}, we deduce that
the following sets are upper semiregularities of $\mathcal{C}(
\mathbb{X}, L(X)).$

\begin{enumerate}

\item $  \mathcal{USR}_1=  \{\mathcal{T} \in \mathcal{C}(
\mathbb{X}, \Phi(X)) : \forall x \in  \mathbb{X}, \,
ind(\mathcal{T}_x) \geq 0 \}$

\item $  \mathcal{USR}_2=  \{\mathcal{T} \in \mathcal{C}(
\mathbb{X}, \Phi(X)): \forall x \in  \mathbb{X}, \,
ind(\mathcal{T}_x) \leq 0  \}$

\item $  \mathcal{USR}_3=  \{\mathcal{T} \in \mathcal{C}(
\mathbb{X}, \Phi(X)) : \forall x \in \mathbb{X}, \,
ind(\mathcal{T}_x) \in m_x \mathbb{Z} \} $

\item $  \mathcal{USR}_4=  \{\mathcal{T} \in \mathcal{C}(
\mathbb{X}, \Phi(X)) : \forall x \in \mathbb{X}, \,
ind(\mathcal{T}_x) = 0 \} $


\item $  \mathcal{USR}_5=  \{\mathcal{T} \in \mathcal{C}(
\mathbb{X}, \Phi_+(X)) : \forall x \in \mathbb{X}, \,\,
ind(\mathcal{T}_x) \leq 0 \}$

\item $  \mathcal{USR}_6=  \{\mathcal{T} \in \mathcal{C}(
\mathbb{X}, \Phi_-(X)) : \forall x \in \mathbb{X}, \,
ind(\mathcal{T}_x) \geq 0 \}$


\end{enumerate}

In fact $\mathcal{USR}_4= \mathcal{USR}_1 \cap \mathcal{USR}_2, $
is the upper semiregularity of Weyl families (see \cite{P49} for
more details about  Weyl families).

\eex

In the case of a single operator, the spectra corresponding to the
upper semiregularities  $ \mathcal{USR}_5$ and $ \mathcal{USR}_6,$
were called respectively the essentiel approximate point spectrum
and the essentiel defect spectrum. See \cite{RAK} and the
references cited there.

{\bf \underline{Open question}:} In  \cite[Theorem 6.5.2]{HA}, it
is proved that every Weyl operator $T$ can be written as $ T=
V+F,$ where $V$ is an invertible operator and $F$ a finite rank
one.

Now let $ \mathcal{T} \in  \mathcal{C}( \mathbb{X}, L(X))$  be  a
Weyl family. Does there exists
  an invertible family $\mathcal{V}  \in \mathcal{C}( \mathbb{X}, L(X)) $  and
 a  finite rank family $\mathcal{F} \in \mathcal{C}( \mathbb{X}, L(X))$
 such that  $ \mathcal{T= V + F}?$

{\small
\noindent Mohammed Berkani,\\
 \noindent Science faculty of Oujda,\\
\noindent University Mohammed I,\\
\noindent Laboratory LAGA, \\
\noindent Morocco\\
\noindent berkanimo@aim.com\\

\end{document}